%% file: kempecolorings.tex
\documentclass{article}
\usepackage{graphicx,color}
\usepackage{amssymb}
\usepackage[geometry]{ifsym}
\usepackage{rotating}

\newtheorem{lemma}{Lemma}
\newtheorem{theorem}{Theorem}

\newtheorem{problem}{Problem}
\newtheorem{conjecture}{Conjecture}
\newcommand{\DONOTTEX}[1]{}

\begin{document}

\title{{\bf Unique colorability and clique minors}}
\author{Matthias Kriesell}

\maketitle

\begin{abstract}
  \setlength{\parindent}{0em}
  \setlength{\parskip}{1.5ex}

  For a graph $G$, let $h(G)$ denote the largest $k$ such that $G$ has $k$ pairwise disjoint pairwise adjacent connected nonempty subgraphs,
  and let $s(G)$ denote the largest $k$ such that $G$ has $k$ pairwise disjoint pairwise adjacent connected subgraphs of size $1$ or $2$.
  {\sc Hadwiger}'s conjecture states that $h(G) \leq \chi(G)$, where $\chi(G)$ is the chromatic number of $G$.
  {\sc Seymour} conjectured $s(G) \geq |V(G)|/2$ for all graphs without antitriangles, i.~e.~ three pairwise nonadjacent vertices.
  Here we concentrate on graphs $G$ with exactly one $\chi(G)$-coloring. We prove generalizations of
  (i) if $\chi(G) \leq 6$ and $G$ has exactly one $\chi(G)$-coloring then $h(G) \geq \chi(G)$, where the proof does {\em not} use  the four-color-theorem, and
  (ii) if $G$ has no antitriangles and $G$ has exactly one $\chi(G)$-coloring then $s(G) \geq |V(G)|/2$.


  
\end{abstract}

\maketitle

\section{Introduction}

All graphs throughout are assumed to be finite, simple, and undirected unless otherwise stated, and for terminology not defined here we refer to 
\cite{BondyMurty2007} or \cite{Diestel2010}.
An {\em anticlique} of a graph $G$ is a set of pairwise nonadjacent vertices of $G$, a {\em coloring} of $G$ is a partition of $V(G)$ 
into anticliques, and $\chi(G)$ denotes the {\em chromatic number} of $G$, that is, the smallest integer $k$ such that $G$ admits a coloring of size $k$.
A {\em clique minor} of $G$ is a set of pairwise disjoint pairwise adjacent connected nonempty subsets of $V(G)$,
where two subsets of $V(G)$ are {\em adjacent} if there exists an edge containing a vertex of each of them, and a subset of $V(G)$ is {\em connected}
if it induces a connected subgraph of $G$. A clique minor is {\em shallow} if all its members have size $1$ or $2$.
Let $h(G)$ be the largest $k$ such that $G$ admits a clique minor of size $k$, and let $s(G)$ denote the largest $k$ such that
$G$ admits a shallow clique minor of size $k$. 

{\sc Hadwiger} conjectured $h(G) \geq \chi(G)$ for all graphs \cite{Hadwiger1943}. So far, this is known for graphs $G$ with $\chi(G) \leq 6$, 
where the statements restricted to $\chi(G)=5$ or to $\chi(G)=6$ are equivalent to the four-color-theorem, respectively \cite{RobertsonSeymourThomas1993}. 
{\sc Seymour} conjectured $s(G) \geq |V(G)|/2$ for all graphs without an {\em antitriangle}, that is, an anticlique of size $3$ (see \cite{Blasiak2007}).
Here, we target these conjectures under the additional assumption of {\em unique colorability}, 
that is, there exists exactly one $\chi(G)$-coloring for the graphs $G$ under consideration.
In the case of {\sc Hadwiger}'s conjecture, we do not get new facts but show that the four-color-theorem is not essential:

\begin{theorem}
  \label{T1}
  If $\chi(G) \leq 6$ and $G$ admits exactly one coloring of size $\chi(G)$ then $h(G) \geq \chi(G)$,
  and the proof of this statement does not rely on the four-color-theorem.
\end{theorem}

Instead of assuming unique colorability, we look at a more general coloring concept:
A coloring ${\frak C}$ of $G$ is a {\em Kempe-coloring} if any two distinct members from ${\frak C}$ induce a connected subgraph of $G$.
It is far from being true that a graph has a Kempe-coloring of any size at all; however, if $G$ has only one coloring ${\frak C}$ of size $k$ then ${\frak C}$ is a
Kempe-coloring, for if, for distinct $A,B$ from ${\frak C}$, $G(A \cup B)$ had more than one component, then we take one, say $H$,
and observe that  ${\frak D}:=({\frak C} \setminus \{A,B\}) \cup \{(A \setminus V(H)) \cup (B \cap V(H)),(B \setminus V(H)) \cup (A \cap V(H))\}$
is a coloring of size $k$ distinct from ${\frak C}$, contradiction. ${\frak D}$ is obtained from ${\frak C}$ by ``exchanging colors'' along the
``Kempe-chain'' $H$, an absolutely classic process in graph coloring theory. The term {\em Kempe-coloring} thus simply indicates that we cannot
apply it as to obtain new colorings of the same size. Now it is clear that Theorem \ref{T1} is a consequence of the following.

\begin{theorem}
  \label{T2}
  If $k \leq 6$ and $G$ admits a Kempe-coloring of size $k$ then $h(G) \geq k$,
  and the proof of this statement does not rely on the four-color-theorem.
\end{theorem}

Whereas I could not relax $k \leq 6$ in Theorem \ref{T1} or Theorem \ref{T2}, the situation improves for {\sc Seymour}'s conjecture.

\begin{theorem}
  \label{T3}
  If $G$ has no antitriangle and exactly one coloring of size $\chi(G)$ then $s(G) \geq \chi(G)$; in particular, $s(G) \geq |V(G)|/2$.
\end{theorem}

It is quite obvious that the colorings of size $\chi(G)$ of a graph $G$ without antitriangles correspond to the maximum matchings of
its complementary graph $\overline{G}$, so that matching theory is naturally involved.
I could not relax the condition of having exactly one coloring of size $\chi(G)$ in Theorem \ref{T3}
to the condition of just having a Kempe-coloring of size $\chi(G)$. However, there is the following, more general ``rooted version'' of Theorem \ref{T3}.
Recall that a {\em transversal} of a set ${\frak S}$ of sets is a set $T$ with $|T \cap D|=1$ for every $D \in {\frak S}$.
If $T$ is a transversal of ${\frak S}$ then we also say that {\em ${\frak S}$ is traversed by $T$}.

\begin{theorem}
  \label{T4}
  Suppose that $G$ has no antitriangles and exactly one coloring ${\frak C}$ of size $\chi(G)$. Then,
  for every transversal $T$ of ${\frak C}$, there exists a shallow clique minor of size $\chi(G)$ traversed by $T$.
\end{theorem}

The paper is organized as follows. In Section \ref{S2} we give a proof of Theorem \ref{T2}, in Section \ref{S3} we prove Theorem \ref{T4};
finally, we discuss some questions and problems in Section \ref{S4}.

\section{Unique colorability and Hadwiger'$\!$s conjecture}
\label{S2}

The main ingredience for the proof of Theorem \ref{T2} is the following result of {\sc Fabila-Monroy} and {\sc Wood} on
``rooted $K_4$-minors'' \cite{FabilaMonroyWood2013}.

\begin{theorem} \cite[Theorem 8,\,(1.)$\leftrightarrow$(3.)]{FabilaMonroyWood2013}
  \label{T5}
  Let $G$ be a graph and $T$ be a set of four vertices of $G$. 
  Then there exists a clique minor of size four traversed by $T$ if and only if
  for any two $a \not= b$ from $T$ there exists a path $P$ in $G$ such that the two vertices from $T \setminus \{a,b\}$ are in the same component of $G-V(P)$.
\end{theorem}
 
The proof of Theorem \ref{T5} in \cite{FabilaMonroyWood2013} is almost self-contained and, in particular, does not use the four-color-theorem.
Only little progress has been made so far on generalizing Theorem \ref{T5} to ``rooted $K_5$-minors'';
it is not clear if there is a set of reasonable ``linkage conditions'' as there at all \cite{Wood2015}.

{\bf Proof of Theorem \ref{T2}.}
We first prove the statement for $k=6$. Let ${\frak C}=\{D_1,\dots,D_6\}$ be a Kempe-coloring of size $6$.
$A:=D_5 \cup D_6$ is connected and has size at least $2$.
There exists a vertex $x \in H$ such that $A_6:=A \setminus \{x\}$ is connected and nonempty
(take a leaf of any spanning tree of $G(A)$). Without loss of generality,
we may assume $x \in D_5$. $A_5:=\{x\}$ is connected and nonempty, too, and adjacent to $A_6$.
Observe that, for $i \not= j$ from $\{1,\dots,6\}$, every vertex from $D_i$ must have at least one neighbor in $D_j$.
It follows that $x$ has a neighbor $b_i$ in $D_i$ for every $i \in \{1,2,3,4\}$.
Let $i \not= j$ from $\{1,2,3,4\}$, and let $i',j'$ be the two indices in $\{1,2,3,4\} \setminus \{i,j\}$.
There is a $b_i,b_j$-path $P$ in $G(D_i \cup D_j)$ and a $b_{i'},b_{j'}$-path $P'$ in $G(D_{i'} \cup D_{j'})$.
Clearly, $P,P'$ are disjoint.
By Theorem \ref{T5}, $G(D_1 \cup D_2 \cup D_3 \cup D_4)$ has a clique minor
$\{A_1,A_2,A_3,A_4\}$ with $b_i$ in $A_i$ for all $i \in \{1,2,3,4\}$. Consequently, $\{A_1,\dots,A_5\}$ is a clique minor,
and, since $b_i \in D_i$ has a neighbor in $D_6 \subseteq A_6$, $\{A_1,\dots,A_6\}$ is a clique minor, too.
This proves $h(G) \geq 6$. --- 

For $k < 6$, let $G^+$ be obtained from $G$ by adding a set $X$ of $6-k$ new vertices and connecting every $x \in X$ to
every other vertex from $V(G) \cup X$ by a new edge. (Later we will refer to the vertices of $X$ as {\em apex vertices}.)
Then ${\frak C}^+:={\frak C} \cup \{\{x\}:\,x \in X\}$ is a Kempe-coloring of $G^+$
of size $6$ so that, by what we have just proved, $G^+$ has a clique minor ${\frak K}$ of size $6$. Since every vertex from $X$
is contained in at most one member of ${\frak K}$, $\{A \in {\frak K}:\,A \cap X=\emptyset\}$ is a clique minor of size at least $6-|X|=k$ of $G$.
Since the $4$-color-theorem is neither used in the proof of Theorem \ref{T5} nor in the preceeding arguments,
this proves Theorem \ref{T2}. \hspace*{\fill}$\Box$

\section{Unique colorability and rooted shallow clique minors}
\label{S3}

For the proof of Theorem \ref{T4}, we need the following result from matching theory by {\sc Kotzig} \cite{Kotzig1959}.
(For a short proof of a more general result let me refer to \cite{Yeo1997}.)
Recall that an edge $e$ is a {\em bridge} of a graph $G$, if $C-e$ is disconnected for some component $C$ of $G$.

\begin{theorem} \cite{Kotzig1959}
  \label{T6}
  If $M$ is the only perfect matching of some nonempty graph $G$ then $M$ contains a bridge of $G$.
\end{theorem}

{\bf Proof of Theorem \ref{T4}.}
We do induction on $|V(G)|$. Let $G$, ${\frak C}$, and $T$ be as in the statement, and set $k:=\chi(G)$.
Observe that ${\frak C}$ is a Kempe-coloring of size $k$, and all members of ${\frak C}$ have order $1$ or $2$.
If some member of ${\frak C}$ consisted of a single vertex, say, $x$, only, then $x \in T$ and $G(D \cup \{x\})$ must be
a star with center $x$ for all $D \in {\frak C} \setminus \{\{x\}\}=:{\frak C}^-$; it follows that $x$ is adjacent to {\em all} other vertices in $V(G)$.
It is easy to see that $\chi(G-x)=k-1$, and that ${\frak C}^-$ is the unique coloring of size $k-1$ of $G-x$.
Clearly, $T^-:=T \setminus \{x\}$ is a transversal of ${\frak C}^-$
so that, by induction, $G-x$ admits a shallow clique minor ${\frak K}^-$ traversed by $T^-$.
But then one readily checks that ${\frak K}:={\frak K}^- \cup \{\{x\}\}$ is a shallow clique minor of $G$ traversed by $T$.

Hence all members of ${\frak C}$ have size $2$, and, thus, they correspond to the edges of a perfect matching $M_0$
of the complementary graph $\overline{G}$ of $G$.
Conversely, if $\overline{G}$ had another perfect matching
$M \not= M_0$, then the edges of $M$ corresponded to the classes of a
coloring of size $k$ of $G$ distinct from ${\frak C}$, a contradiction.
Therefore, $M_0$ is the only perfect matching of $\overline{G}$.

We now construct a descending sequence of subgraphs $H_0,\dots,H_\ell$ of $\overline{G}$ with unique perfect matchings as follows.
Set $H_0:=\overline{G}$, suppose that $H_0,\dots,H_\ell$ have been constructed, and let $M_\ell$ be the unique perfect matching of $H_\ell$.
By Theorem \ref{T6}, $M_\ell$ contains a bridge $e_\ell$ of $H_\ell$. Let $Q_\ell$ be a component of $H_\ell-e_\ell$ containing exactly one
vertex of $e_\ell$ such that $|V(Q_\ell)| \leq |V(H_\ell)|/2$, and let $x_\ell,y_\ell$ be the vertices of $e_\ell$, where $x_\ell \in V(Q_\ell)$.
Observe that $xy \in E(G)$ for all $(x,y) \in (V(Q_\ell) \times (V(H_\ell) \setminus V(Q_\ell))) \setminus \{(x_\ell,y_\ell)\}$.
It is easy to see that $R_\ell:=\{e \in M_\ell:\,V(e) \subseteq V(Q_\ell)\}$ is the unique perfect matching of $Q_\ell-x_\ell$ (which is possibly empty),
and that $M_\ell \setminus (R_\ell \cup \{e_\ell\})$ is the unique perfect matching of $H_\ell \setminus (V(Q_\ell) \cup V(e_\ell))$;
by the size condition to $Q_\ell$, the latter matching contains at least as many edges as $R_\ell$, and
we choose a subset $S_\ell$ of size exactly $|R_\ell|$ from it. Set $T_\ell:=H_\ell(\bigcup \{ V(e): \, e \in S_\ell\} \cup \{y_\ell\})$
(so $S_\ell$ is the unique perfect matching of $T_\ell-y_\ell$).  Set $H_{\ell+1}:=H_\ell-(V(Q_\ell) \cup V(T_\ell))$;
then $M_{\ell+1}:=M_\ell-(R_\ell \cup S_\ell \cup \{e_\ell\})$ is the unique perfect matching of $H_{\ell+1}$.
We iterate until $H_{\ell+1}$ is empty.

\begin{figure}
  \begin{center} 
    \scalebox{0.5}{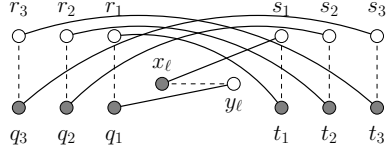} 
    \caption{\label{F3} Constructing ${\frak K}_\ell$ if $|R_\ell|=|S_\ell|=z=3$ and $x_\ell \in T$. Dashed lines are the
    edges of $R_\ell \cup S_\ell \cup \{e_\ell=x_\ell y_\ell\}$, solid lines correspond to the members of ${\frak K}_\ell$.
    Vertices from $T$ are displayed grey.}
  \end{center}
\end{figure}

If $R_\ell \not=\emptyset$ then we can construct a set ${\frak K}_\ell$ of cliques of size $2$ of $G$
traversed by $T$ such that each contains both a member of $V(Q_\ell)$ and a member of $V(T_\ell)$:
Let $q_1r_1,\dots,q_zr_z$ be the edges of $R_\ell$ and $s_1t_1,\dots,s_zt_z$ be those of $S_\ell$,
where $q_1,\dots,q_z$ and $t_1,\dots,t_z$ are from $T$; if $x_\ell \in T$ then set
${\frak K}_\ell:=\{\{x_\ell,s_1\}$, $\{y_\ell,q_1\}$, $\{r_1,t_1\}\} \cup \{\{q_j,s_j\},\{r_j,t_j\}:\,j \in \{2,\dots,z\}\}$,
and if, otherwise, $y_\ell \in T$ then set
${\frak K}_\ell:=\{\{x_\ell,t_1\},\{y_\ell,r_1\},\{q_1,s_1\}\} \cup \{\{q_j,s_j\},\{r_j,t_j\}:\,j \in \{2,\dots,z\}\}$
(see Figure \ref{F3} for an example of the first kind).
Otherwise, $S_\ell=R_\ell=\emptyset$, so that $Q_\ell$ consists of $x_\ell$ only; if $x_\ell \in T$ then we set
${\frak K}_\ell:=\{\{x_\ell\}\}$, and in the other case we call the index $\ell$ {\em special} and set ${\frak K}_{\ell}:=\emptyset$.
Observe that each ${\frak K}_\ell$ is a shallow clique minor of $G(V(Q_\ell) \cup V(T_\ell))$ traversed by $T$.

We first apply induction to $G^-:=\overline{H_1}$ ($=$ $G(V(H_1))$, with the unique coloring
of size $\chi(G^-)$ corresponding to the edges of $M_1$, traversed by $T \cap V(G^-)$), and find
a shallow clique minor ${\frak K}^-$ of $G^-$ traversed by $T \cap V(G^-)$.

Suppose that $0$ is not a special index. Since every member of ${\frak K}_0$ contains a vertex
from $V(Q_0)$ and all of these are adjacent to all of $V(G^-)$, ${\frak K}_0 \cup {\frak K}^-$ is a shallow clique minor of $G$ traversed by $T$,
and we are done.
Hence we may assume that $0$ is a special index. If $\{y_0\}$ is adjacent to all members from ${\frak K}^-$ in $G$
then ${\frak K}^- \cup \{\{y_0\}\}$ is a shallow clique minor of $G$ traversed by $T$. Hence we may assume
that $y_0$ is nonadjacent to at least one vertex $u \in T \setminus \{y_0\}$. There exists a vertex $v \in V(G)$ such that $uv \in M_0$;
clearly, $v$ is  not in $T$ as $u$ is in $T$, and $v$ is adjacent to $y_0$ as $G$ contains no antitriangles.
There exists a unique $i>0$ such that $v \in V(Q_i) \cup V(T_i)$, and we choose $v$ (and $u$) such that $i$ is as large as possible. 
If $V(Q_i) \cup V(T_i)=\{x_i,y_i\}=V(e_i)$ then let $w$ be the vertex 
from $V(e_i)$ distinct from $v$;  otherwise, there exists a set $A \in {\frak K}_i$ such that $v \in A$, and we let
$w$ be the vertex in $A \setminus \{v\}$.
Observe that, in both cases, $w \in T$ follows.
By induction, applied to $G^-:=\overline{H_{i+1}}$, there exists a shallow clique minor ${\frak K}^-$ of size $\chi(G^-)=|V(G^-)|/2$ traversed by $T$.

\begin{figure}
  \begin{center}
    \scalebox{0.31}{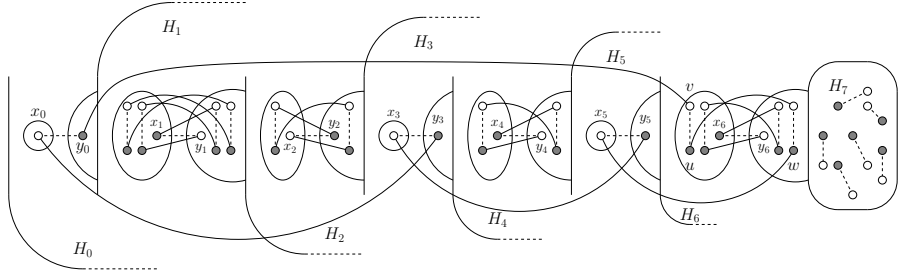} 
    \caption{\label{F1}  Example for a decomposition like in the proof of Theorem \ref{T4}, with $i=6$ and special indices $0,3,5$.
      Dashed edges are from the matching $M_0$ (among them the ``horizontal'' edges $e_0,\dots,e_6$),
      solid edges are from $G$ and constitute part of the shallow clique minor,
      grey vertices are from the transversal $T$, white ones are from $V(G) \setminus T$.
      The subgraphs $Q_0,\dots,Q_6$ are encircled, and their vertices are adjacent to everything to their right, respectively.}
  \end{center}
\end{figure}

Now let $0:=s_0,s_1,\dots,s_d$ be the special indices smaller than $i$ in increasing order (see Figure \ref{F1} for an example). 
We let $A_d:=\{w,x_{s_d}\}$, and $A_j:=\{y_{s_{j+1}},x_{s_j}\}$ for $j \in \{0,\dots,d-1\}$, and $A^+:=\{y_0,v\}$.
Since $x_{s_d}$ and $w \in H_i$ are in different components of $H_{s_d}-e_{s_d}$ and $w \not= y_{s_d}$, $A_d$ is a clique.
Since $x_{s_j}$ and $y_{s_{j+1}}$ are in different components of $H_{s_j}-e_{s_j}$ and $y_{s_{j+1}} \not= {y_{s_j}}$, $A_j$ is a clique for $j \in \{0,\dots,d-1\}$.
By construction, $A^+$ is a clique. All of $A_0,\dots,A_d,A^+$ are pairwise disjoint, disjoint from any clique in ${\frak K}^-$, and
disjoint from any clique in ${\frak K}_i$ with the only exception that $A_d$ meets the unique clique of ${\frak K}_i$ containing $w$
(which is either $\{v,w\}$ or $\{w\}$) if $i$ is not special.
Observe that ${\frak K}^- \cup \bigcup_{j=0}^i {\frak K}_j$ has size $|V(G)|/2-(d+1)$ if $i$ is not special, and size $|V(G)|/2-(d+2)$ if $i$ is special.
Therefore,  ${\frak K}:={\frak K}^- \cup \bigcup_{j=0}^{i-1} {\frak K}_j \cup ({\frak K}_i \setminus \{\{v,w\},\{w\}\}) \cup \{A_0,A_1,\dots,A_d,A^+\}$
consists of exactly $|V(G)|/2$ many disjoint cliques of size $1$ or $2$. Moreover, ${\frak K}$ is traversed by $T$,
as ${\frak K}^-$, each ${\frak K}_j$, and $\{A_0,\dots,A_d,A^+\}$ are traversed by $T$, respectively.

We accomplish the proof by showing that ${\frak K}$ is a clique minor. To this end, it suffices to prove for all $\ell \geq 0$ that if
$A \in {\frak K}$ contains some vertex from $T \cap (V(Q_\ell) \cup V(T_\ell))$ and
$B \in {\frak K} \setminus \{A\}$ contains some vertex from $T \cap V(H_\ell)$ then $A,B$ are adjacent.
This is true for $\ell>i$, as $A,B$ are members of the clique minor ${\frak K}^-$ in this case.

Consider $\ell \leq i$ and $A \not=B$ from ${\frak K}$ with
$A \cap T \cap (V(Q_\ell) \cup V(T_\ell)) \not= \emptyset$ and $B \cap T \cap V(H_\ell) \not=\emptyset$.
Suppose, to the contrary, that $A,B$ are not adjacent in $G$.
If $A$ is contained in ${\frak K}_\ell$ then it contains a vertex from $V(Q_\ell)$,
which is adjacent in $G$ to all vertices from $V(H_\ell) \setminus (V(Q_\ell) \cup \{y_\ell\})$,
implying that $B \subseteq V(Q_\ell) \cup \{y_\ell\}$ unless $B=A_j$ for some $j \in \{0,\dots,d\}$.
If $B=A_j$ then $x_{s_j} \in B \setminus V(H_\ell)$ and $s_j<\ell$; since $x_{s_j}$ is adjacent to all of $V(H_{s_j}) \setminus \{x_{s_j},y_{s_j}\}$,
it is adjacent to all of $A \subseteq V(H_\ell) \subseteq V(H_{s_j})$, contradiction. Hence $B \subseteq V(Q_\ell) \cup \{y_\ell\}$.
If $A=\{x_\ell\}$ then $R_\ell=S_\ell=\emptyset$ and ${\frak K}_\ell=\{\{x_\ell\}\}$,
so that $B=\{y_\ell\}$ is not traversed by $T$, contradiction. Thus, $A \in {\frak K}_\ell$ has size $2$ and contains a vertex from $V(T_\ell)$, which is
adjacent in $G$ to all vertices from $V(Q_\ell) \setminus \{x_\ell\}$, so that $B \subseteq \{x_\ell,y_\ell\}$.
Since $x_\ell,y_\ell$ are nonadjacent in $G$, $B=\{x_\ell\}$ or $B=\{y_\ell\}$ follows. We infer $B=\{x_j\}=V(Q_\ell)$
(as all of $A_0,\dots,A_d,A^+$ have size $2$), so that $A$ does not contain a member of $V(Q_\ell) \setminus B$, contradiction.
Therefore, $A$ must be among $A_0,\dots,A_d,A^+$. 

If $A=A_j$ for some $j \in \{0,\dots,d\}$ then $A$ contains $x_{s_j}$ from $V(G) \setminus T$, and $A$ contains $y_{s_{j+1}}$ from $T$ if $j<d$ or $w$
from $T$ if $j=d$. In both cases, $s_j<\ell$, and $x_{s_j}$ from $A$ is adjacent to all of $V(H_{s_j}) \setminus \{x_{s_j},y_{s_j}\}$,
and, thus, to all of $B \subseteq V(H_\ell) \subseteq V(H_{s_j})$, contradiction.
If, finally, $A=A^+=\{y_0,v\}$ then $\ell=0$. By choice of $i$, $y_0$ is nonadjacent to all vertices from $V(H_{i+1}) \setminus T$ and, hence,
adjacent to all of $V(H_{i+1}) \cap T$, implying that $B \not\in {\frak K}^-$.
Since every vertex from $V(Q_j)$, $j<i$, is adjacent to $v$, we know that $A=A^+$ is adjacent to all $B$ from $\bigcup_{j=0}^{i-1} {\frak K}_j$
and to all of $A_0,\dots,A_d$. It follows that $B \in {\frak K}_i \setminus \{\{w\},\{v,w\}\}$. But then $B$ has a vertex in $V(Q_i)$ and
a vertex in $V(T_i)$; if $v \in V(T_i)$ then it is adjacent to all of $V(Q_i)$, and if $v \in V(Q_i)$ then it is adjacent to all of $V(T_i)$,
so that in either case, $v$ is adjacent to $B$.
\hspace*{\fill}$\Box$

\section{Open problems}
\label{S4}

So far, I could neither generalize Theorem \ref{T2} to the case $k>6$, nor Theorem \ref{T4} under the weaker assumption that
there is a Kempe-coloring of order $\chi(G)$ instead of a unique coloring of size $\chi(G)$, not even if the conclusion asks for
just {\em any} clique minor of size $\chi(G)$ instead of a shallow one.
The following, being a common generalization of these two projects, is, therefore, perhaps a little bit too optimistic.

\begin{conjecture}
  \label{C1}
  Suppose that $G$ has a Kempe-coloring ${\frak C}$ of size $k$.
  Then, for every transveral $T$ of ${\frak C}$, there exists 
  a clique minor traversed by $T$.
\end{conjecture}

For $k=4$, this follows from Theorem \ref{T5}, as in the proof of Theorem \ref{T2} in Section \ref{S2},
and the statement inherits to smaller values of $k$
by augmenting $G,{\frak C}$ to $G^+,{\frak C}^+$ as above (add only $4-k$ instead of $6-k$ apex vertices and augment $T$ to $T^+:=T \cup X$).
It suffices to prove the Conjecture for the case that $G$ is (edge-) minimal with the property that ${\frak C}$ is a Kempe-coloring of size $k$.
This is equivalent to saying that any two distinct members of ${\frak C}$ induce a tree, which is in turn
equivalent to saying that --- independent from the actual Kempe-coloring of size $k$ --- the graph has
exactly $|E(G)|=(k-1) \cdot |V(G)| - {k \choose 2}$ edges.
Since, for fixed $|V(G)|$, the latter is an increasing function in $k$ as long as $k \leq |V(G)|$,
it follows easily that the minimality condition implies that $G$ has no Kempe-coloring of size larger than $k$ (proofs are left to the reader).

\begin{figure}
  \begin{center} 
    \scalebox{0.7}{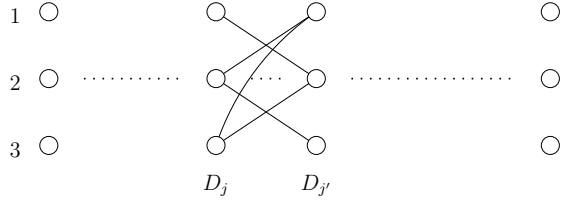} 
    \caption{\label{F2} For $j<j'$, $D_j \cup D_{j'}$ induces a path on six vertices.}
  \end{center}
\end{figure}

In general, an affirmative answer to Conjecture \ref{C1} could give us a clique minor of much larger size than
we actually need to verify {\sc Hadwiger}'s conjecture for the graphs under consideration,
because the difference of the size of a (largest) Kempe-coloring (if any) and $\chi(G)$ may be arbitrarily large:
Consider the graph $G$ on $V(G):=\{1,2,3\} \times \{1,\dots,k\}$, $k \geq 3$, where $(i,j) \not= (i',j')$ are connected
if and only if $i \not= i'$ and $j \not= j'$ and ($i=2 \vee i'=2 \vee (i>i' \wedge j<j'))$.
Consider $D_j:=\{(i,j):\,i \in \{1,2,3\}\}$. For $j<j'$, $D_j \cup D_{j'}$ induces a path on six vertices (see Figure \ref{F2}), which implies
that ${\frak C}:=\{D_j:\,j \in \{1,\dots,k\}\}$ is a Kempe-coloring of size $k$ and $G$ is minimal with that property.
(By the remark above, ${\frak C}$ would even be a largest Kempe-coloring of $G$.)
On the other hand, $\{\{(i,j):\, j \in \{1,\dots,k\}\}:\,i \in \{1,2,3\}\}$ is a coloring of size $3$, and the triangle
formed by $(3,1)$, $(1,2)$, $(2,3)$ shows $\chi(G)=3$.
Adding $\ell$ apex vertices to $G$ as in the proof of Theorem \ref{T2} yields a graph
$G^+$ with a Kempe-coloring of size $k+\ell$ and chromatic number $3+\ell$.

The proof of Theorem \ref{T4} (and also the much easier proof of the unrouted version, where special indices can be avoided) depends
heavily on {\sc Kotzig}'s Theorem, Theorem \ref{T5}, and the latter does not generalize to the situation of a Kempe-coloring because
there are infinitely many graphs $G$ without antitriangles admitting a Kempe-coloring of size $\chi(G)=|V(G)|/2$
whose complementary graph is bridgeless and, in fact, $(|V(G)|+2)/4$-connected: Let $D$ be a $k$-connected tournament on $2k+1$ vertices.
We construct $G$ from $D$ by introducing a pair $x^+,x^-$ for every vertex $x$ from $D$ and connect $x^+$ and $y^-$
if $x=y$ or there is a directed edge from $x$ to $y$ in $D$. It is easy to see that $G$ is a $(k+1)$-connected bipartite graph on $4k+2$ vertices.
The perfect matching $M:=\{x^+x^-:x \in V(D)\}$ of $G$ corresponds to a coloring ${\frak C}$ of size $2k+1$ in $\overline{G}$,
and, since distinct members of $M$ are connected by only one edge in $G$, any pair of members of ${\frak C}$ induces a connected graph
(in fact: a path of length $3$) in $\overline{G}$. Therefore, ${\frak C}$ is a Kempe-coloring of $\overline{G}$ (even a largest one).
--- An answer to the following question could lead further:

\begin{problem}
  \label{P1}
  Let $G$ be a triangle free graph with a perfect matching $M$ such that any two edges of
  $M$ are connected by at most one edge (out of the four possible edges).
  Does there exists an $A \subseteq V(G)$ with $|A|=|V(G)|/2$ such that every $4$-cycle of $G$ contains
  two nonadjacent vertices from $A$ or two nonadjacent vertices from $V(G) \setminus A$?
\end{problem}

The answer is ,,yes'' if, for example, $G$ is bipartite (like the graphs above obtained from tournaments) 
or if the girth of $G$ is larger than $4$.
(In fact, I do not know a single triangle-free graph on an even number of vertices where the answer is ,,no'', but I should doubt that there aren't any.)
The point is, of course, that if $G$ is a graph without antitriangles providing a Kempe-coloring of size $k=|V(G)|/2$ then $\overline{G}$
meets the conditions in Problem \ref{P1} (where the matching $M$ of $\overline{G}$ corresponds to the Kempe-coloring of $G$).
Suppose that we get a set $A$ as in Problem \ref{P1}.
We may assume that $\overline{G}$ is $k=|V(\overline{G})|/2$-connected (see below), so that, 
by {\sc Hall}'s or {\sc Menger}'s Theorem (see \cite{Diestel2010} or \cite{BondyMurty2007}),
there exists a perfect matching $N$ from $A$ to $V(G) \setminus A$ in $\overline{G}$. If, for distinct $e,f$ from $N$,
all four potential edges between $e,f$ were absent in $\overline{G}$, then $V(e) \cup V(f)$ would induce a $4$-cycle in $G$ without two nonadjacent vertices 
from $A$ or two nonadjacent vertices from $V(G) \setminus A$, contradiction.
So the edges of $M$ correspond to a shallow clique minor of $\overline{G}$.

Concerning the connectivity issue, let us prove the following Lemmas.

\begin{lemma}
  \label{L1}
  Suppose that $G$ has a Kempe-coloring ${\frak C}$ and let $T \subseteq V(G)$ be a separator of $G$.
  Then $T$ contains an element of all but at most one member of ${\frak C}$.
\end{lemma}

{\bf Proof.}
Suppose that there exists a $D \in {\frak C}$ such that $D \cap T=\emptyset$. Then there exists a component $C$ of $G-T$ such
that $V(C) \cap D \not= \emptyset$. Let $C' \not= C$ be another component of $G-T$. If $V(C') \cap D \not= \emptyset$
then $T$ contains a member of every $D' \in {\frak C} \setminus \{D\}$, since $G(D \cup D')$ is connected.
Otherwise, $D \subseteq V(C)$, and there exist $x \in V(C')$ and $B \in {\frak C} \setminus \{D\}$ with $x \in B$. 
Since $G(B \cup D)$ is connected, $T$ contains
a member of $B$. For any $D' \in {\frak C} \setminus \{B,D\}$, $x$ has a neighbor $y \in D'$. 
Since $y \in T \cup V(C')$ and $G(D \cup D')$ is connected, $T$ contains a member of $D'$, too. 
\hspace*{\fill}$\Box$

\begin{lemma}
  \label{L2}
  Suppose that $G$ has no antitriangles and has a Kempe-coloring of size $k$.
  Then $G$ is $k$-connected or admits a shallow clique minor of size $k$.
\end{lemma}

{\bf Proof.}
Suppose that $G$ is a minimal counterexample to the statement of the Lemma, and let ${\frak C}$ be a Kempe-coloring of $G$ of size $k$.
Since $G$ is not a clique of size $k$, $G$ had a separator $T$ of size $k-1$ by the previous Lemma.
Among all separators of size $k-1$, we choose $T$ and a component $C$ of $G-T$ such that 
$V(C)$ is minimal (with respect to either set inclusion or size).

If some member of ${\frak C}$ consisted of a single vertex $x$ then, as usual, $x$ would be adjacent to all other vertices,
$T \setminus \{x\}$ would be a separator of $G-x$ of size $k-2$,  and ${\frak C} \setminus \{\{x\}\}$ would be a Kempe-coloring of size $k-1$
of $G-x$. By minimality of $G$, $G-x$ had a shallow clique minor ${\frak K}^-$ of size $k-1$, so that
${\frak K}^- \cup \{\{x\}\}$ would be a shallow clique minor of size $k$ of $G$, contradiction.

Hence all members of ${\frak C}$ have size $2$. There exists a $D \in {\frak C}$ with $D \cap T = \emptyset$.
By Lemma \ref{L1}, $|B \cap T|=1$ for all $B \in {\frak C} \setminus \{D\}$. Since $G-T$ has no antitriangles,
$G-T$ has exactly one further component $C'$ distinct from $C$, and $V(C')$ and $V(C)$ induce cliques. 
We classify ${\frak C} \setminus \{D\}$ by 
${\frak B}:=\{B \in {\frak C} \setminus \{D\}:\, B \cap V(C) \not= \emptyset\}$ and 
${\frak B}':=\{B \in {\frak C} \setminus \{D\}:\, B \cap V(C') \not = \emptyset\}$.
If $D \subseteq V(C)$ then there exists a $B \in {\frak B}'$ as $V(C') \not= \emptyset$, but the vertex in $B \cap V(C')$
cannot have a neighbor in $D$, contradicting the fact that $G(B \cup D)$ is connected.
It follows that there exist unique vertices $d \in D \cap V(C)$ and $d' \in D \cap V(C')$.

Let $P':=\bigcup {\frak B}'$. If there was no matching from $P' \cap T$ into $V(C')$
then, by {\sc Hall}'s Theorem (see \cite{BondyMurty2007} or \cite{Diestel2010}), $|N_G(X) \cap V(C')|<|X|$ for some $X \subseteq P' \cap T$, so that
$S:=(T \setminus X) \cup (N_G(X) \cap V(C'))$ contained the neighborhood of $Y:=V(C') \setminus N_G(X)$.
Since $Y$ is a nonempty (proper) subset of $V(C')$ and $|S| \leq k-2$, $S$ is a separator of $G$ of size less than $k-1$, contradiction.
Therefore, we find a matching $M'$ from $P' \cap T$ into $V(C')$. 

Let $P:=\bigcup {\frak B}$. If there was no matching from $P \cap T$ into $P \cap V(C)$
then, by {\sc Hall}'s Theorem,
$|N_G(X) \cap P \cap V(C)|<|X|$ for some $X \subseteq P \cap T$, so that
$S:=(T \setminus X) \cup (N_G(X) \cap P \cap V(C)) \cup \{d\}$ contained the neighborhood of $Y:=P \cap V(C) \setminus N_G(X)$.
Since $Y$ is a nonempty proper subset of $V(C)$ and $|S| \leq k-1$, $S$ is a separator of $G$ of size $k-1$ and $Y$ contains a component
of $G-S$ whose vertex set is properly contained in $V(C)$, contradicting the choice of $C$.
Therefore, we find a matching $M$ from $P \cap T$ into $P \cap V(C)$. (It may be that $P=M=\emptyset$.)

The edges of $M$ are pairwise adjacent, since $V(C)$ is a clique, and 
the edges of $M'$ are pairwise adjacent, since $V(C')$ is a clique.
Now let $e \in M$ and $e' \in M'$. Then there exists a $B \in {\frak B}$ and $x \in V(e) \cap T \cap B$
and a $B' \in {\frak B}'$ and $x' \in V(e') \cap T \cap B'$. Since $G(B \cup B')$ is connected and lacks the edge connecting
the two vertices in $(B \cup B') \setminus T$, we know that $xx' \in E(G)$.
Therefore, $e,e'$ are adjacent.
Now $d$ is not an endvertex of any member of $M \cup M'$, but it is adjacent to all of $M$ as $V(C)$ is a clique
and to all of $M'$ as, for every $B' \in {\frak B}'$,
$d$ must have a neighbor in $B'$ which can only be the vertex of $B' \cap T$.

It follows that $\{V(e):\, e \in M \cup M'\} \cup \{\{d\}\}$ is a shallow clique minor of size $k$ in $G$.
\hspace*{\fill}$\Box$

According to the above arguments, verifying Conjecture \ref{C1} for some value of $k$ would verify the statement of Theorem \ref{T2} for $k+2$;
therefore, I think it would be already interesting to verify Conjecture \ref{C1} for the smallest open case of $k=5$.
In general, graphs on $n$ vertices with a Kempe-coloring of size $5$ must have at least $4n-10$ edges;
by a classic result of {\sc Thomassen} \cite{Thomassen1974}, they admit not only a minor but even a subdivision of $K_5$ where
one could, moreover, prescribe a single vertex in the interior of one of the subdivision paths. By a result of {\sc Mader} \cite{Mader1998}
(answering a question of {\sc Dirac}), $3n-5$ edges suffice to guarantee a clique minor of size $5$
(implying that there is such a minor in any graph $G$ which has an edge such that $G-e$ has a Kempe-coloring of size $4$).
Both results indicate that there is considerable freedom in choosing a clique minor of size $5$ provided that $G$ has a Kempe-coloring of size $5$.

\DONOTTEX{ 
Another related question would ask for a suitable generalization of {\sc Kotzig}'s Theorem, Theorem \ref{T5}, as follows:
Given a connected graph with a unique smallest partition ${\frak M}$ of $V(G)$ into cliques, is there a $C \in {\frak M}$ such that $G-E(G(C))$
is disconnected? 
}

{\bf Author's Address.}

{\sc Matthias Kriesell} \\
Technische Universit\"at Ilmenau \\
Weimarer Stra{\ss}e 25 \\
D--98693 Ilmenau \\
Germany

\end{document}

%% file: kempe3.pdf_t
\begin{picture}(0,0)%
\includegraphics{kempe3.pdf}%
\end{picture}%
\setlength{\unitlength}{3947sp}%
\begingroup\makeatletter\ifx\SetFigFont\undefined%
\gdef\SetFigFont#1#2#3#4#5{%
  \reset@font\fontsize{#1}{#2pt}%
  \fontfamily{#3}\fontseries{#4}\fontshape{#5}%
  \selectfont}%
\fi\endgroup%
\begin{picture}(4678,1927)(2312,-7157)
\put(6901,-5461){\makebox(0,0)[b]{\smash{{\SetFigFont{17}{20.4}{\rmdefault}{\mddefault}{\updefault}{\color[rgb]{0,0,0}$s_3$}%
}}}}
\put(5701,-5461){\makebox(0,0)[b]{\smash{{\SetFigFont{17}{20.4}{\rmdefault}{\mddefault}{\updefault}{\color[rgb]{0,0,0}$s_1$}%
}}}}
\put(6301,-5461){\makebox(0,0)[b]{\smash{{\SetFigFont{17}{20.4}{\rmdefault}{\mddefault}{\updefault}{\color[rgb]{0,0,0}$s_2$}%
}}}}
\put(5701,-7061){\makebox(0,0)[b]{\smash{{\SetFigFont{17}{20.4}{\rmdefault}{\mddefault}{\updefault}{\color[rgb]{0,0,0}$t_1$}%
}}}}
\put(6301,-7061){\makebox(0,0)[b]{\smash{{\SetFigFont{17}{20.4}{\rmdefault}{\mddefault}{\updefault}{\color[rgb]{0,0,0}$t_2$}%
}}}}
\put(3601,-7061){\makebox(0,0)[b]{\smash{{\SetFigFont{17}{20.4}{\rmdefault}{\mddefault}{\updefault}{\color[rgb]{0,0,0}$q_1$}%
}}}}
\put(2401,-5461){\makebox(0,0)[b]{\smash{{\SetFigFont{17}{20.4}{\rmdefault}{\mddefault}{\updefault}{\color[rgb]{0,0,0}$r_3$}%
}}}}
\put(3001,-5461){\makebox(0,0)[b]{\smash{{\SetFigFont{17}{20.4}{\rmdefault}{\mddefault}{\updefault}{\color[rgb]{0,0,0}$r_2$}%
}}}}
\put(3601,-5461){\makebox(0,0)[b]{\smash{{\SetFigFont{17}{20.4}{\rmdefault}{\mddefault}{\updefault}{\color[rgb]{0,0,0}$r_1$}%
}}}}
\put(2401,-7061){\makebox(0,0)[b]{\smash{{\SetFigFont{17}{20.4}{\rmdefault}{\mddefault}{\updefault}{\color[rgb]{0,0,0}$q_3$}%
}}}}
\put(3001,-7061){\makebox(0,0)[b]{\smash{{\SetFigFont{17}{20.4}{\rmdefault}{\mddefault}{\updefault}{\color[rgb]{0,0,0}$q_2$}%
}}}}
\put(6901,-7061){\makebox(0,0)[b]{\smash{{\SetFigFont{17}{20.4}{\rmdefault}{\mddefault}{\updefault}{\color[rgb]{0,0,0}$t_3$}%
}}}}
\put(5101,-6661){\makebox(0,0)[b]{\smash{{\SetFigFont{17}{20.4}{\rmdefault}{\mddefault}{\updefault}{\color[rgb]{0,0,0}$y_\ell$}%
}}}}
\put(4201,-6161){\makebox(0,0)[b]{\smash{{\SetFigFont{17}{20.4}{\rmdefault}{\mddefault}{\updefault}{\color[rgb]{0,0,0}$x_\ell$}%
}}}}
\end{picture}%

%% file: kempe2.pdf_t
\begin{picture}(0,0)%
\includegraphics{kempe2.pdf}%
\end{picture}%
\setlength{\unitlength}{3947sp}%
\begingroup\makeatletter\ifx\SetFigFont\undefined%
\gdef\SetFigFont#1#2#3#4#5{%
  \reset@font\fontsize{#1}{#2pt}%
  \fontfamily{#3}\fontseries{#4}\fontshape{#5}%
  \selectfont}%
\fi\endgroup%
\begin{picture}(18024,5424)(889,-9073)
\put(9601,-6061){\makebox(0,0)[b]{\smash{{\SetFigFont{17}{20.4}{\rmdefault}{\mddefault}{\updefault}{\color[rgb]{0,0,0}$y_3$}%
}}}}
\put(13801,-6061){\makebox(0,0)[b]{\smash{{\SetFigFont{17}{20.4}{\rmdefault}{\mddefault}{\updefault}{\color[rgb]{0,0,0}$y_5$}%
}}}}
\put(4201,-4261){\makebox(0,0)[b]{\smash{{\SetFigFont{20}{24.0}{\rmdefault}{\mddefault}{\updefault}{\color[rgb]{0,0,0}$H_1$}%
}}}}
\put(7501,-8461){\makebox(0,0)[b]{\smash{{\SetFigFont{20}{24.0}{\rmdefault}{\mddefault}{\updefault}{\color[rgb]{0,0,0}$H_2$}%
}}}}
\put(9301,-4561){\makebox(0,0)[b]{\smash{{\SetFigFont{20}{24.0}{\rmdefault}{\mddefault}{\updefault}{\color[rgb]{0,0,0}$H_3$}%
}}}}
\put(13201,-4861){\makebox(0,0)[b]{\smash{{\SetFigFont{20}{24.0}{\rmdefault}{\mddefault}{\updefault}{\color[rgb]{0,0,0}$H_5$}%
}}}}
\put(17701,-5461){\makebox(0,0)[b]{\smash{{\SetFigFont{20}{24.0}{\rmdefault}{\mddefault}{\updefault}{\color[rgb]{0,0,0}$H_7$}%
}}}}
\put(7501,-6161){\makebox(0,0)[b]{\smash{{\SetFigFont{17}{20.4}{\rmdefault}{\mddefault}{\updefault}{\color[rgb]{0,0,0}$y_2$}%
}}}}
\put(8701,-5961){\makebox(0,0)[b]{\smash{{\SetFigFont{17}{20.4}{\rmdefault}{\mddefault}{\updefault}{\color[rgb]{0,0,0}$x_3$}%
}}}}
\put(12901,-5961){\makebox(0,0)[b]{\smash{{\SetFigFont{17}{20.4}{\rmdefault}{\mddefault}{\updefault}{\color[rgb]{0,0,0}$x_5$}%
}}}}
\put(2401,-8761){\makebox(0,0)[b]{\smash{{\SetFigFont{20}{24.0}{\rmdefault}{\mddefault}{\updefault}{\color[rgb]{0,0,0}$H_0$}%
}}}}
\put(2401,-6661){\makebox(0,0)[b]{\smash{{\SetFigFont{20}{24.0}{\rmdefault}{\mddefault}{\updefault}{\color[rgb]{0,0,0}$y_0$}%
}}}}
\put(14701,-8061){\makebox(0,0)[b]{\smash{{\SetFigFont{20}{24.0}{\rmdefault}{\mddefault}{\updefault}{\color[rgb]{0,0,0}$H_6$}%
}}}}
\put(10801,-8161){\makebox(0,0)[b]{\smash{{\SetFigFont{20}{24.0}{\rmdefault}{\mddefault}{\updefault}{\color[rgb]{0,0,0}$H_4$}%
}}}}
\put(1501,-5961){\makebox(0,0)[b]{\smash{{\SetFigFont{20}{24.0}{\rmdefault}{\mddefault}{\updefault}{\color[rgb]{0,0,0}$x_0$}%
}}}}
\put(14701,-5461){\makebox(0,0)[b]{\smash{{\SetFigFont{20}{24.0}{\rmdefault}{\mddefault}{\updefault}{\color[rgb]{0,0,0}$v$}%
}}}}
\put(16801,-7061){\makebox(0,0)[b]{\smash{{\SetFigFont{20}{24.0}{\rmdefault}{\mddefault}{\updefault}{\color[rgb]{0,0,0}$w$}%
}}}}
\put(14701,-7061){\makebox(0,0)[b]{\smash{{\SetFigFont{20}{24.0}{\rmdefault}{\mddefault}{\updefault}{\color[rgb]{0,0,0}$u$}%
}}}}
\put(3901,-6161){\makebox(0,0)[b]{\smash{{\SetFigFont{17}{20.4}{\rmdefault}{\mddefault}{\updefault}{\color[rgb]{0,0,0}$x_1$}%
}}}}
\put(6601,-6661){\makebox(0,0)[b]{\smash{{\SetFigFont{17}{20.4}{\rmdefault}{\mddefault}{\updefault}{\color[rgb]{0,0,0}$x_2$}%
}}}}
\put(4801,-6661){\makebox(0,0)[b]{\smash{{\SetFigFont{17}{20.4}{\rmdefault}{\mddefault}{\updefault}{\color[rgb]{0,0,0}$y_1$}%
}}}}
\put(15301,-6161){\makebox(0,0)[b]{\smash{{\SetFigFont{17}{20.4}{\rmdefault}{\mddefault}{\updefault}{\color[rgb]{0,0,0}$x_6$}%
}}}}
\put(16201,-6661){\makebox(0,0)[b]{\smash{{\SetFigFont{17}{20.4}{\rmdefault}{\mddefault}{\updefault}{\color[rgb]{0,0,0}$y_6$}%
}}}}
\put(10801,-6161){\makebox(0,0)[b]{\smash{{\SetFigFont{17}{20.4}{\rmdefault}{\mddefault}{\updefault}{\color[rgb]{0,0,0}$x_4$}%
}}}}
\put(11701,-6661){\makebox(0,0)[b]{\smash{{\SetFigFont{17}{20.4}{\rmdefault}{\mddefault}{\updefault}{\color[rgb]{0,0,0}$y_4$}%
}}}}
\end{picture}%

%% file: kempe1.pdf_t
\begin{picture}(0,0)%
\includegraphics{kempe1.pdf}%
\end{picture}%
\setlength{\unitlength}{3947sp}%
\begingroup\makeatletter\ifx\SetFigFont\undefined%
\gdef\SetFigFont#1#2#3#4#5{%
  \reset@font\fontsize{#1}{#2pt}%
  \fontfamily{#3}\fontseries{#4}\fontshape{#5}%
  \selectfont}%
\fi\endgroup%
\begin{picture}(4904,1767)(1786,-2939)
\put(1801,-2561){\makebox(0,0)[b]{\smash{{\SetFigFont{12}{14.4}{\rmdefault}{\mddefault}{\updefault}{\color[rgb]{0,0,0}$3$}%
}}}}
\put(4501,-2861){\makebox(0,0)[b]{\smash{{\SetFigFont{12}{14.4}{\rmdefault}{\mddefault}{\updefault}{\color[rgb]{0,0,0}$D_{j'}$}%
}}}}
\put(1801,-1361){\makebox(0,0)[b]{\smash{{\SetFigFont{12}{14.4}{\rmdefault}{\mddefault}{\updefault}{\color[rgb]{0,0,0}$1$}%
}}}}
\put(1801,-1961){\makebox(0,0)[b]{\smash{{\SetFigFont{12}{14.4}{\rmdefault}{\mddefault}{\updefault}{\color[rgb]{0,0,0}$2$}%
}}}}
\put(3601,-2861){\makebox(0,0)[b]{\smash{{\SetFigFont{12}{14.4}{\rmdefault}{\mddefault}{\updefault}{\color[rgb]{0,0,0}$D_j$}%
}}}}
\end{picture}%